\documentclass[12pt]{article}
\usepackage{amssymb}
\usepackage{graphics}
\usepackage{epsfig}
\usepackage{harvard}
\citationstyle{dcu}

\makeatletter

\@addtoreset{equation}{section}

\makeatother

\newtheorem{thm}[equation]{Theorem}

\newtheorem{propos}[equation]{Proposition}

\def\reff#1{(\ref{#1})}
\def\sobre#1#2{\lower 1ex \hbox{ $#1 \atop #2 $ } }

\begin{document}

\def\E{{\mathbb E}}
\def\P{{\mathbb P}}
\def\R{{\mathbb R}}
\def\Z{{\mathbb Z}}
\def\V{{\mathbb V}}
\def\N{{\mathbb N}}
\def\NN{{\bf N}}
\def\X{{\cal X}}
\def\supp{{\rm Supp}\,}
\def\XX{{\bf X}}
\def\Y{{\bf Y}}
\def\G{{\cal G}}
\def\T{{\cal T}}
\def\cC{{\cal C}}
\def\C{{\bf C}}
\def\D{{\bf D}}
\def\U{{\bf U}}
\def\K{{\bf K}}
\def\H{{\bf H}}
\def\n{{\bf n}}
\def\m{{\bf m}}
\def\b{{\bf b}}
\def\g{{\bf g}}
\def\sqr{\vcenter{
         \hrule height.1mm
         \hbox{\vrule width.1mm height2.2mm\kern2.18mm\vrule width.1mm}
         \hrule height.1mm}}                  
\def\square{\ifmmode\sqr\else{$\sqr$}\fi}
\def\one{{\bf 1}\hskip-.5mm}
\def\liml{\lim_{L\to\infty}}
\def\given{\ \vert \ }
\def\ze{{\zeta}}
\def\be{{\beta}}
\def\de{{\delta}}
\def\la{{\lambda}}
\def\ga{{\gamma}}
\def\th{{\theta}}
\def\proof{\noindent{\bf Proof. }}
\def\rate{{e^{- \beta|\ga|}}}
\def\A{{\bf A}}
\def\B{{\bf B}}
\def\C{{\bf C}}
\def\D{{\bf D}}
\def\MM{{\bf m}}
\def\lnt{{\Lambda^N}}
\def\S{{\mathcal{S}}}
\def\basis{{\rm Basis}\,}
\def\life{{\rm Life}\,}
\def\birth{{\rm Birth}\,}
\def\death{{\rm Death}\,}
\def\flag{{\rm Flag}\,}
\def\color{{\rm Color}\,}
\def\type{{\rm Type}\,}

\title{Improved bounds for perfect simulation of a continuous
  one-dimensional loss network} \author{Nancy L. Garcia, {\it
  UNICAMP} \thanks{Corresponding author: IMECC/UNICAMP, Caixa Postal
  6065, 13.081-970 - Campinas SP BRAZIL, {\tt 
  nancy@ime.unicamp.br}}  \\ Nevena Mari\'c, {\it USP} \thanks{
  IME/USP, Caixa Postal 66281, 05311-970 - S\~{a}o Paulo SP BRAZIL,
  {\tt nevena@ime.usp.br}}} 
\date{\today}

\maketitle

\abstract Perfect simulation of an one-dimensional loss network on
$\R$ with length distribution $\pi$ and cable capacity $C$ is
performed using the clan of ancestors method. Domination of the clan of
ancestors by a branching process with longer memory improves the
sufficient conditions for the perfect scheme to be
applicable. 
\endabstract

\section{Introduction} \label{sec:intro}

Kelly (1991) introduced a continuous unbounded loss network described
as follows. Imagine that users are arranged along an infinitely long
cable and that a call between two points on the cable $s_{1}$, $s_{2}
\in {\mathbb R}$ involves just that section of the cable between $s_{1}$
and $s_{2}$. Past any point along its length the cable has the
capacity to carry simultaneously up to $C$ calls: a call attempt
between $s_{1}$, $s_{2} \in {\mathbb R}$, $s_{1} < s_{2}$, is lost if
past any point of the interval $[s_{1}, s_{2}]$ the cable is already
carrying $C$ calls. Suppose that calls are attempted at points in
${\mathbb R}$ following a homogeneous Poisson process with rate $\lambda$.
Assume that the section of the cable demanded by a call has distribution
$\pi$ with finite mean $\rho_1$ and the duration of a call has
exponential distribution with mean one. Assume that the location of a
call, the cable section needed and its duration are independent.  Let
$m(s,t)$ be the number of calls in progress past point $s$ on the
cable at time $t$. Kelly (1991) conjectured that $((m(s,t), s \in {\mathbb
  R}), t \ge 0)$ has a unique invariant measure, given by a stationary
$M/G/\infty$ queue (Markov arrivals, 
general service time and infinite servers) conditioned to have at most $C$
clients at all times. Ferrari and Garcia (1998) used a continuous
(non-oriented) percolation argument to prove the above conjecture
whenever $\pi$ has finite third moment and 
the arrival rate $\lambda$ is sufficiently small. Fern\'{a}ndez,
Ferrari and Garcia (2002) using an oriented percolation argument
improved this bound to 
\begin{equation}
\label{eq:t54a}
\lambda_c^{FFG} = \frac{1}{(\rho_2 + \rho_1 +1)}
\end{equation}
where $\rho_1$ and $\rho_2$ are the first and second moment of
distribution $\pi$ respectively. This argument is based on a graphical
representation of the birth and death process and it is the basis for
the perfect simulation scheme ``\emph{Backward-Forward Algorithm}'',
described in Fern\'{a}ndez, Ferrari and Garcia (2002). This algorithm
involves the ``thinning'' of a marked Poisson process ---the
\emph{free process}--- which dominates the birth-and-death process,
and it involves a time-backward and a time-forward sweep. The initial
stage of the construction is done \emph{toward the past}, starting
with a finite window and retrospectively looking to \emph{ancestors},
namely to those births in the past that could have (had) an influence
on the current birth.  The construction of the \emph{clan of
  ancestors} constitutes the time-backward sweep of the algorithm.
Once this clan is completely constructed, the algorithm proceeds in a
time-forward fashion ``cleaning up'' successive generations according
to appropriate penalization schemes. The relation ``being ancestor
of'' induces a backward in time \emph{contact/oriented percolation}
process.  The algorithm is applicable as long as this oriented
percolation process is sub-critical. Garcia and Mari\'c (2003) using
the Perron-Frobenius theory for sub-criticality of branching process
obtained a new bound given by
\begin{equation}
\label{eq:t54b}
\lambda_c^*=\frac{1}{(\sqrt{\rho_2} + \rho_1)}.
\end{equation}
However, studying the characteristics of the clan of ancestors through
simulation in Section \ref{sec:simul} it is clear that the domination
by the branching process is not sharp. That is, the number of
ancestors is much smaller than the total number of the population in
the branching process and the clan of ancestors can be finite even
though the branching is supercritical. By studying the perfect
simulation algorithm and constructing a 2-generation branching
process, that is a Markov process with order 2, it is possible to
improve bound \reff{eq:t54b}  to
\begin{equation}
\lambda_c^{**} =    \frac{4}{3\rho_1  +
\sqrt{\rho_1^2 + 8\rho_2 }}
\end{equation}
and it stands that  $\lambda_c^{**} \ge (4/3)
\lambda_c^*$ .

\section{Loss networks} 

Let $\G$ be a
family of intervals of the line $\gamma$ ($\gamma=(x,x+u), x,u\in
\R$), which will be named \emph{calls}, and consider a state space
${\cal S}=\{\xi\in\N^\G\,:\,$ $\xi(\ga)\neq 0$ only for a countable
set of $\ga\in\G\}$.

Loss networks are a particular class of spatial birth-and-death
processes. The evolution of these processes in time
   are given either by the {\em birth} of a new call to be added to the
   actual configuration or by the {\em death} of an existing call that
   will be eliminated from the actual configuration. Moreover, they
   have the Markovian property in time that, the probability of a
   change depends only on the actual configuration of the
   system. A loss network $\eta_t$ is defined by a marked
Poisson process where births are controlled by a {\em birth rate},  a 
   non-negative measurable function $\mathbf{b}(\gamma,\eta)$ such that
   $$\int_B \mathbf{b}(\gamma,\eta)d\gamma<\infty$$
   for each $B$,
   bounded Borel set, and for all $\eta \in {\cal S}$.  The births are
   regulated by the exclusion principle, depending on the capacity
   ($C$) of the network. The marks include a life-time exponentially
   distributed with mean one and the length of the call which is
   distributed according to a distribution $\pi$ with finite first and
   second moments $\rho_1$ and $\rho_2$ respectively. The \emph{death
     rate} also a non-negative measurable function $\mathbf{d}:
   \mathbb{R}^d\times {\cal S} \rightarrow [0,\infty).$ In this work,
   we are going to assume that the death rate is always equal to one.
   The generator of the process is given by
\begin{equation}
Af(\eta)=\int
{(f(\eta+\delta_\gamma)-f(\eta))b(\gamma,\eta)d\gamma} +
\int(f(\eta-\delta_\gamma)-f(\eta))\eta(d\gamma)
\end{equation}
where $\eta \in  \{0,1\}^{\mathcal{B}(\mathbb{R})}$. The death rate
$1$ is included in the second expression.

For a process $\alpha_t$ with rate densities which are
independent of the 
actual configuration there exists $\omega: \G \rightarrow
[0,\infty)$ such that
\begin{equation}
  \label{eq:l40}
  \mathbf{b}(\gamma,\alpha)=\omega(\gamma)
\end{equation} 
we call this process a \emph{free process}.  Such a process is just a
space-time marked Poisson process.  It exists and is ergodic whichever
the choice of $w$. In the particular case where
$\omega(\gamma)=\lambda$ the invariant measure is the
$\lambda$-homogeneous Poisson process. For the one-dimensional loss
networks (described in Section \ref{sec:intro}) the birth rate is uniformly 
bounded and it can be decomposed as 
\begin{equation} \label{decomp}
\mathbf{b}((x,x+u),\eta)=\lambda ~\pi(u)~M((x,x+u),\eta)
\end{equation}
where, $0\leq M(\gamma,\eta)\leq 1$. The first factor represents a
basic birth-rate density due to an ``internal'' Poissonian clock and
the last factor acts as an unnormalized probability for the individual to
be actually born once the internal clock has rang.  The birth
is hindered or reinforced according to the configuration $\eta$.

In this case, for capacity $C=1$,
\begin{equation}
M(\gamma,\eta)=\prod_{\theta: \eta(\theta)\neq
0}(1-I(\gamma,\theta))
\end{equation}
\begin{equation} \label{indic}
I(\gamma,\theta)= \left\{ \begin{array}{cc}
                          1 &  \gamma\cap \theta \neq \emptyset\\
                          0 & \mbox{otherwise}
                          \end{array}
                          \right.
\end{equation}
where $\gamma, \theta$ are of the form  $(x,x+u)$. For $C>1$,
the expression is less simple
\[
M((x,x+u),\eta)= \left\{
                          \begin{array}{cc}
                          1 &  \mbox{otherwise} \\
                          0 & \mbox{ there exists } y\in (x,x+u)~~ and\\
                            &  \theta_1,\ldots,\theta_C \mbox{ such that }
\eta(\theta_i)=1\\
                            &   \mbox{ and } y\in \theta_i~ \mbox{ for all
 } i=1,\ldots,C.
                          \end{array}
                          \right. \]

\section{Graphical construction for the loss networks}\label{graphical}

Let $\cal{N}$ = $\{$ ($\xi_1,T_1),$ ($\xi_2,T_2)$,\ldots$\}$ be a homogeneous
Poisson Process with rate $\lambda$ in $\mathbb{R}\times
[0,\infty)$,  $S_1,S_2,\ldots$ be i.i.d.  random variables
exponentially distributed with mean one and  $U_1, U_2,\ldots$
 be i.i.d. random variables with common distribution $\pi$. Assume the
 family of variables  $\{S_1, S_2, \ldots\}$, $\{U_1, U_2, \ldots\}$
 and the Poisson process are all independent. Consider the random
 rectangles  
\begin{center}
$R_i=\{(x,y); \xi_{i}\leq x\leq \xi_{i}+U_i, T_{i}\leq y\leq
T_{i}+S_i\}$.
\end{center}

Then $\{R_i, i\geq 1\}= \{(\xi_{i},T_i)+D_i, i \geq1\}$ is a Boolean
model in  $\mathbb R^2$ 
where $D_i = [0,U_i]\times [0,S_i]$ and represents the free process of
calls. Boolean models have the property that the number of sets $C \in
\cal C$ that cover a fixed point $x \in \R^d$ is a Poisson random
variable with mean $\lambda \E(\mbox{vol}(S))$. For more details about
coverage processes see Hall (1988).

Now, for each rectangle $R_i$ we associate an independent mark $Z_i
\sim U(0,1)$, and each marked rectangle we identify with the marked
point  $(\xi_i, T_i, S_i, U_i, Z_i)$.
We recognize in the marked point process 
 $\mathbf{R}=\{(\xi_i,T_i, S_i, U_i, Z_i),~i=1,2,\dots \}$ a graphical
 representation of the birth and death process with constant birth
 rate  $\lambda$, and constant death rate, equal to 1. We call this
 free  process  $\alpha$ and  $Z_i$ will serve as a flag of allowed
 births. Calling $R = (\xi,\tau,s,u,z)$, we use the notation 
\begin{equation}
\basis(R)= (\xi, \xi+u), ~~  \birth(R)= \tau, ~~  \life(R)=
[\tau,\tau+s], ~~  \flag(R)= z.
\end{equation}
We also define, for two rectangles $R$ and $R'$, 
\begin{center}
$R'\nsim R$, if $R'\cap R \neq \emptyset$\\ $R'\sim R$, otherwise.
\end{center}

We need a series of definitions:
\begin{itemize}
\item For an arbitrary point $(x,t)\in \mathbb{R}^2$ define the
  collection of all rectangles in $\mathbf R$ that contain this point 
\begin{equation}
\A_1^{(x,t)}=\{R\in \mathbf{R}| ~x\in \basis(R), t\in \life(R)\}
\end{equation}
\item For each rectangle $R$ define its ancestor set 
\begin{equation}
\A_1^R=\{R'\in \mathbf{R}|~\birth(R')\leq \birth(R),~ R'\nsim R\}
\end{equation}
\item Define recursively the generations ($n > 1$) of the above sets
  that is, the $n$th generation of ancestors:
\begin{eqnarray}
\A_n^{(x,t)}=\{R''| R''\in \A_1^{R'}\mbox{ for some }R'\in
\A_{n-1}^{(x,t)} \} \\ \A_n^{R}=\{R''| R''\in \A_1^{R'} \mbox{ for some
}R'\in \A_{n-1}^{R} \} 
\end{eqnarray}
We say that there is \emph{backward oriented percolation} if there
exists one point   $(x,t)$ such that
$\A_n^{(x,t)}\neq \emptyset$ for all $n$, that is, if there exists
one point with an infinite number of ancestors. Call \emph{clan of
  ancestors} of $(x,t)$ the union of all its ancestors:
\begin{equation}
\A^{(x,t)}=\bigcup_{n\geq 1} \A_n^{(x,t)}
\end{equation}
and $\bf{R}\rm[0,t]=\{R\in \bf{R}\rm|~\birth(R)\in [0,t]\}$.
\end{itemize}

The existence of the process in infinite volume for any time interval
is guaranteed as long as the process do not explode, that is, no
rectangle has an infinite number of ancestors in a finite time. The
following theorem is proved in Fern\'{a}ndez, Ferrari and Garcia (2001).
\begin{thm}
If $\A^{(x,t)}\cap \bf{R}\rm[0,t]$ is finite with probability one, for
any  $x\in \mathbb R$ and $t\geq 0$, then for all 
$\Lambda\subseteq \mathbb R$ the loss network process defined in
$\Lambda$  is well-defined and has at least one
invariant measure $\mu^\Lambda$.
\end{thm}

For the existence of the process in infinite time, it is needed that
the clan of ancestors of all rectangles are finite with probability
one, that is, there is no backward oriented percolation. In order to
construct the invariant measure for stationary Markov processes it is
usual to construct the process beginning at $-\infty$ with an
arbitrary configuration and look at the process at time $0$. If the
configuration at time $0$ does not depend on the initial configuration
then we have a sample of invariant measure. The graphical construction
described above allow us to construct the process $\eta_t$
by a thinning of the free process  $\alpha_t$ for all $t\in
\mathbb{R}$. Moreover, the same argument shows that the distribution
of $\eta_0$ does not depend on the initial configuration. The next
theorem summarizes the results about the process, see  Fern\'andez et
al. (2001, 2002). 
\begin{thm}
  If with probability one there is no backward oriented percolation
  in $\mathbf{R}$, then the loss network process can be constructed in
  $(-\infty,\infty)$ in such a way that the marginal distribution of
  $\eta_t$ is invariant. Moreover, this distribution is unique and the
  velocity of convergence is exponential.
\end{thm}

One way of determining the lack of percolation is the domination
through a branching process. Establishing sub-criticality conditions
for the branching process we obtain sufficient conditions for lack of
percolation. Looking backward, the ancestors will be the branches. The
time of the death will be the birth time for the branching process.
The clan of ancestors in itself is not a branching process because the
lack of independence.

Let $R$ be a rectangle with basis $\gamma = (x,x+u)$ with length 
$u$, born at time  $0$. Define $\tilde{b}_n^{u} (v)$ as the number of
rectangles in the $n$th generation of ancestors of  $R$ having basis
with length  $v$:
\begin{equation}\label{procram}
\tilde{b}_n^{u}(v) = |\{ R^{'} \in \A_n^R |~ |\basis(R^{'})|= v\} |.
\end{equation}

The process $\tilde{b}_n$ is not a Galton-Watson process but it can be
dominated by one (call it $b_n$) as described by Fern\'andez et al. (2001),
where each call length represents a type and it has as offspring
distribution the same one as $\tilde{b}_1$. The number of types can be
finite, countable or uncountable depending upon the distribution
$\pi$.

For the one-dimensional loss network the 
offspring distribution of $b_n$ is Poisson distributed  with mean  
\begin{equation}\label{media} 
m(u, v)= \lambda ~ \pi(v)~(u+ v)
\end{equation}
where $m(u, v)$ is the mean number of children type $v$ for parents
type $u$.  In this case, Garcia and Mari\'c (2003) used the
Perron--Frobenius theory for sub-criticality of branching process
obtained a new bound given by
\begin{equation}
\lambda(\sqrt{\rho_2} + \rho_1) <1.
\end{equation}

In Section \ref{sec:ramificacao} we are going to improve this bound
dominating the clan of ancestors by a branching process with order 2,
that is, the reproducing mechanism is governed not only by the parents
but also by the grandparents.

\section {Backward-Forward Algorithm (BFA) applied to loss networks}

\indent
The BFA was introduced by  Fern\'andez, Ferrari e Garcia (2002) to
perfect simulate from spatial point processes which are absolutely
continuous with respect to a Poisson point process and that are invariant
measures of spatial birth and death processes. 

The algorithm does involve the ``thinning'' of a marked Poisson process
---the \emph{free process}--- which dominates the birth-and-death
process, and it involves a time-backward and a time-forward sweep.
But these procedures are performed in a form quite different from
previous algorithms.  The initial stage of our construction is done
\emph{toward the past}, starting with a finite window and
retrospectively looking to \emph{ancestors}, namely to those births in
the past that could have (had) an influence on the current birth.  The
construction of the \emph{clan of ancestors} constitutes the
time-backward sweep of the algorithm.  Once this clan is completely
constructed, the algorithm proceeds in a time-forward fashion
``cleaning up'' successive generations according to appropriate
penalization schemes.  

The relation ``being ancestor of'' induces a backward in time
\emph{contact/oriented percolation} process.  The algorithm is
applicable as long as this oriented percolation process is
sub-critical.

To simplify the implementation of BFA to the loss network process we
 are going to assume that $\pi$ have compact support. This assumption
 is not necessary and can be removed with a little modification on the
 generation of the free process. Define $H=\inf\{y~|~\pi((0,y))=1\}$.

\subsection{Construction of the clan of ancestors of a finite
 window  $\Lambda=[a,b]\subset \mathbb R$} \label{sec:clan}

We are interested in sampling  a finite window $\Lambda=[a,b]$ of the
  equilibrium  measure in \emph{infinite-volume}. 

\begin{itemize}
\item[C1.] Generate the free process $\alpha_0 =
  \{\xi_1^0,\ldots,\xi_m^0 \}$; a homogeneous Poisson process with rate
  $\lambda$ in the interval $[a-H,b]$.

$s_L^0=a$;~~$s_R^0=b$.

\item[C2.] Generate $U_1^0,\ldots,U_{m}^0$ i.i.d. random variables with
  common distribution $\pi$ and let  $\eta = \emptyset$.
\newline For each  $i$ from 1 to $m$
\begin{equation}
  \mbox{ if } (\xi_i^0,\xi_i^0+U_i^0)\cap [a,b]\neq \emptyset \mbox{ then }
 \eta = \eta \cup (\xi_i^0,\xi_i^0+U_i^0)
\end{equation}
We are simply generating rectangles with basis intersecting
$[a,b]$. We have $n_0=|\eta|\leq m$ basis.

\item[C3.] Generate $S_1^0,\ldots,S_{n_0}^0$ i.i.d. exponential random
  variables with mean one and construct the rectangles
\begin{equation}
\mathbf{R}_0\mathrm = \{(\xi_i^0,\xi_i^0+U_i^0)\times [-S_i^0,0]; i =
1,\ldots,n\}. 
\end{equation}
Consider now the following subset of  $\mathbb{R}\times (-\infty,0]$
\begin{equation}
\Lambda_0=\bigcup_{i=1}^{n} (\xi_i^0-H,\xi_i^0+U_i^0)\times [-S_i^0,0]
\end{equation}

\item[C4.] $k=1$;~$\Delta= \Lambda_0$;

\item[C5.]\label{skok}
$s_L^{k}= \min(s_L^{k-1},\min_{i \leq n_{k-1}}(\xi_i^{k-1} -H))$\\
 $s_R^{k}=\max(s_R^{k-1},\max_{i \leq n_{k-1}}(\xi_i^{k-1} +U_i^{k-1}))$

\item[C6.] Generate a $\lambda$-homogeneous Poisson process  
$\{(\xi_1^k,\tau_1^k),\ldots,(\xi_{n_k}^k,\tau_{n_k}^k)\}$ on
$\Delta \cup [s_L^{k},s_L^{k-1}) \cup (s_R^{k-1}, s_R^{k}]$.

\item[C7.] Generate  $U_1^k,\ldots,U_{n_k}^k$ i.i.d. random variables with
  distribution $\pi$ and  $S_1^k,\ldots,S_{n_k}^k$ i.i.d. exponential random
  variables with mean one and construct the rectangles
\begin{equation}
\mathbf{R}_k\mathrm = \{(\xi_i^k,\xi_i^k+U_i^k)\times
[\tau_i^k-S_i^k,\tau_i^k]; i = 1,\ldots,n_k\}.
\end{equation}
Consider
\begin{equation}
\Lambda_k=\bigcup_{i=1}^{n_k} (\xi_i^k-H,\xi_i^k+U_i^k)\times
[\tau_i^k-S_i^k,0].
\end{equation}

\item[C8.]
\begin{itemize}
\item if $n_k = 0$ then construct the clan of ancestors of $\eta$
\begin{equation}
 \A^{\eta}:=\bigcup_{i=0}^{k-1}\mathbf{R}_i
\end{equation}
and STOP.
\item otherwise, do 
\newline $\Delta = \Lambda_k\backslash \Lambda_{k-1}$;
\newline k=k+1;
\newline return to C5;
\end{itemize}

\end{itemize}

\noindent {\it Remark.} At Step C6., it is necessary to consider only
rectangles satisfying
$$\tau_i^k-S_i^k > \min_{j=1,\dots,n_{k-1}}(\tau_i^{k-1}-S_i^{k-1}).$$
This restriction does not affect the generation of the clan of
ancestors but reduces drastically the computational cost.

\newblock
 We finish performing the  \emph{BACKWARD} step of the algorithm: the
 construction of the clan of ancestors. The \emph{FORWARD} step
 corresponds to move to the beginning of the clan of ancestors and
 decide which rectangles are going to be kept and which ones are going
 to be erased.  Once these clans are perfectly simulated, it is only
necessary to apply the \emph{deterministic} ``cleaning procedure'', based on
the capacity $C$ of the network , to obtain a perfect sample of the
 interacting process. In this case, if a point $(x,t)$ belongs to more
 than $C$ rectangles, keep the $C$ rectangles born first and erase the others.

\subsection{The cleaning algorithm}
Call \bf{T} \rm the set of rectangles to be tested and \bf{K} \rm the
set of kept rectangles. 
\begin{itemize}
\item[L1.] \bf{K}\rm=$\emptyset$; \bf{T}\rm= $\A^\eta$;
\item[L2.] If $\bf{T}\rm=\emptyset$ go to L4.
\newline otherwise, order \bf{T }\rm  by birth time. Let
$R_1$ be the first rectangle following such ordering.
\newline $\bf{K}\rm=\bf{K}\rm\cup R_1$;
$\bf{T}\rm=\bf{T}\rm\backslash R_1$
\item[L3.] Depending upon $C$
\begin{enumerate}
\item If $C=1$; For all $R\in \bf{T}\rm$ such that $R\nsim R_1$,
$\bf{T}\rm=\bf{T}\rm\backslash R$.
\newline return to  L2.
\item Se $C>1$;
\newline \qquad \emph{for i=1 to $|\bf{T}\rm|-C$}
 \newline $R_i\in \bf{T}\rm$, if $R_i\nsim R_1$ call
\emph{Area}=$R_i\cap R_1$, $C(Area)$=2 and
$\bf{K}\rm=\bf{K}\rm\cup R_i$; $\bf{T}\rm=\bf{T}\rm\backslash
R_i$;
\newline \qquad ~~\emph{for j=1 to $|\bf{T}\rm|$}
\newline
\qquad if $R_j\cap \emph{Area}\neq \emptyset$ take
$C(Area)=C(Area)+1$, if $C(Area)>C$ then
$\bf{T}\rm=\bf{T}\rm\backslash R_j$
\newline return to L2.
\end{enumerate}
\item [L4.] Take $\bf{K}\rm ^\eta = \bf{K}\rm$ and STOP.
\end{itemize}

\vskip5mm

Obtaining $\bf{K}\rm ^\eta$, we define
\begin{eqnarray}
\eta^{\ast}(\gamma)= \left\{ \begin{array}{cc}
                          1 &  \eta(\gamma)=1 \mbox{ and } \exists R
                          \in \mathbf{K}\rm ^\eta \mbox{ such that }
                          \basis(R)=\gamma \\ 
                          0 & \mbox{otherwise}
                          \end{array}
                          \right.
\end{eqnarray}
Theorem 3.18 of Fern\'andez et al. (2002) guarantees that
$\eta^{\ast}$ is a perfect sample from the invariant measure of the
loss network described above.

\pagebreak

\subsection{Simulation results}

In this section we present some of the simulation results for several
values of $\lambda$, $C$ ( network capacity) and window $\Lambda$. The
distribution $\mathbf{\pi}$ is taken to be $\mathbf{U(0,1)}$. In this
case, Garcia and Mari\'{c} (2003) obtained that $\lambda<0.9282$ is a
sufficient condition for the simulation. The programs were written in
MATLAB 5.0. For easiness of reading the results are presented in two
steps: the clan of ancestors and the cleaning result. The basis of the
rectangles kept at time $t=0$ constitutes the perfect sample.

\begin{figure}[h]
\centering
\psfig{file = 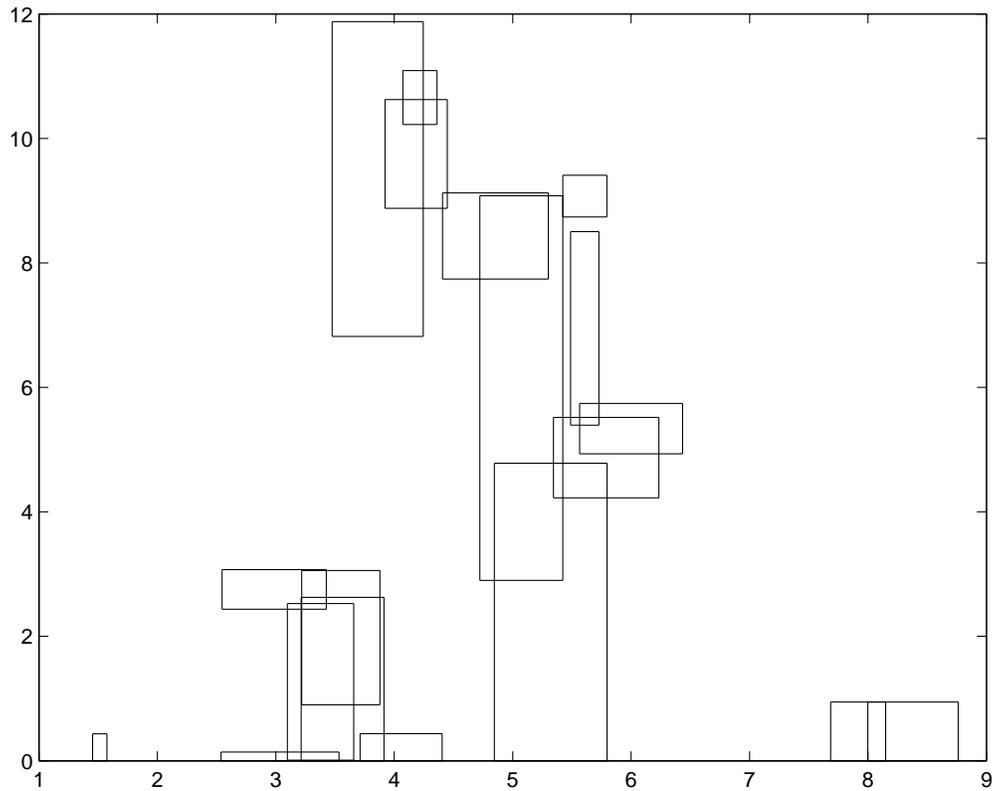,scale=.8}
 \label{c07}
\caption{Clan of ancestors for  $U(0,1)$, $\lambda=0.7$, $\Lambda=[0,10]$.}
\end{figure}

\begin{figure}[p]\centering
\psfig{file = 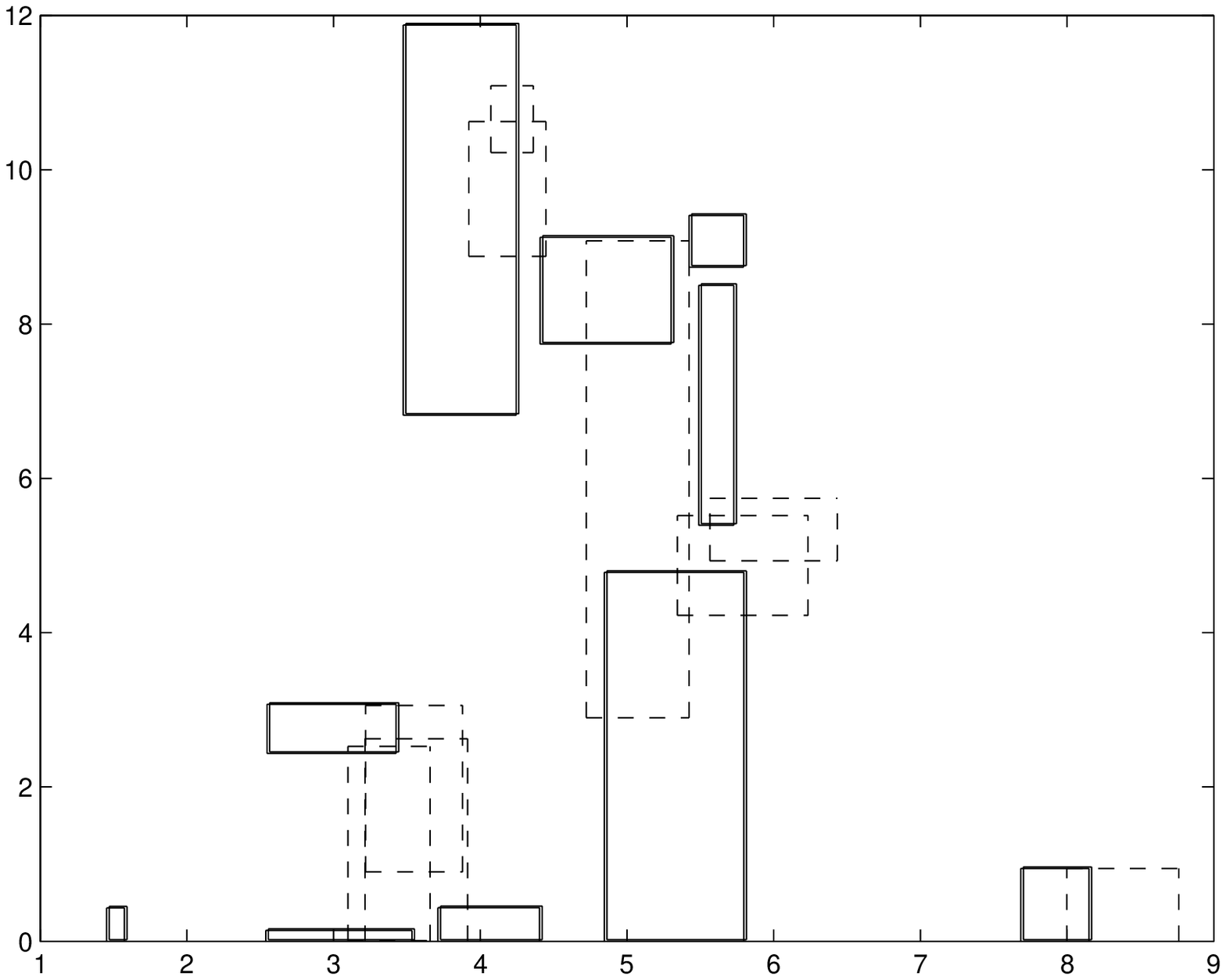,scale=.6}
\psfig{file = 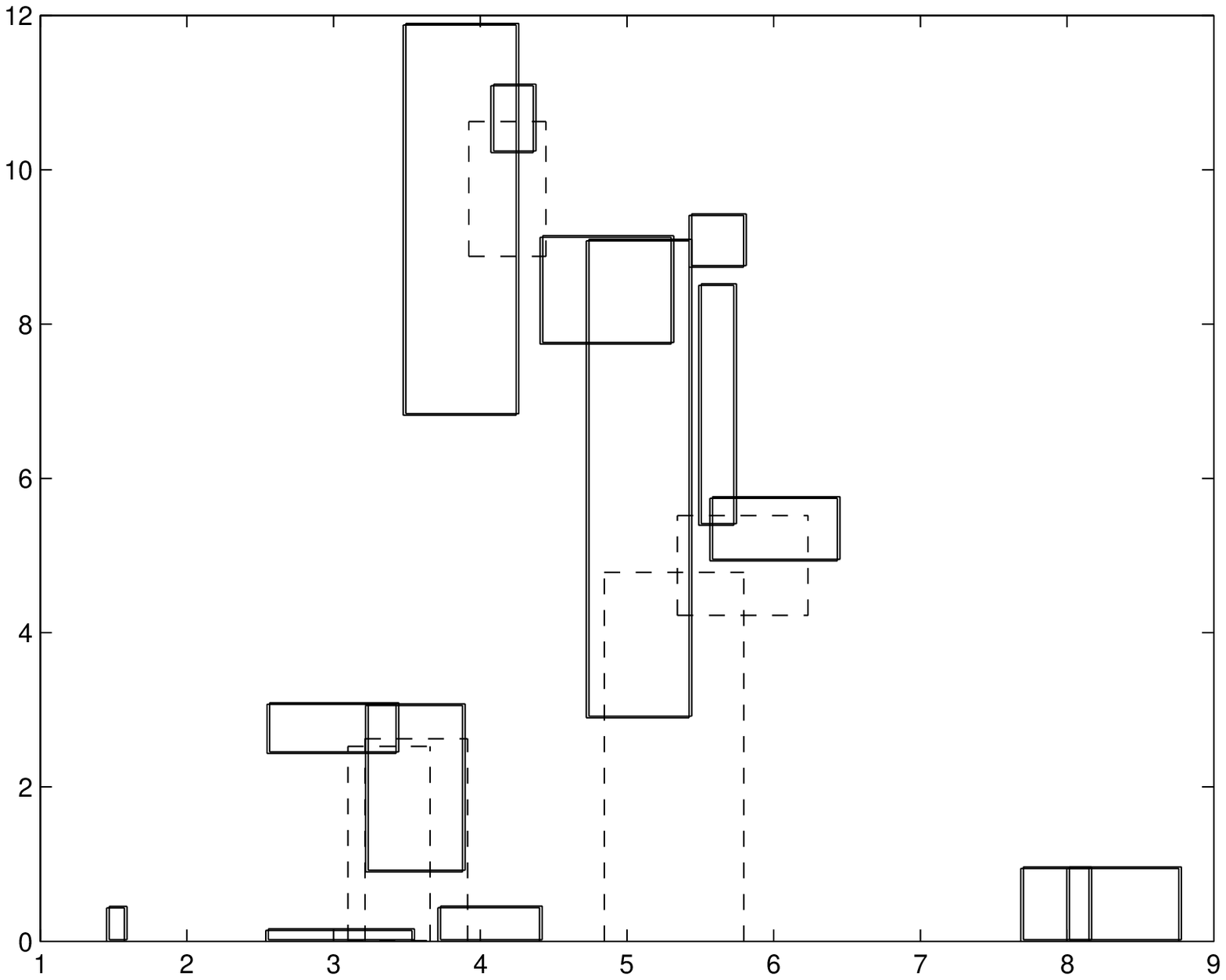,scale=.6}
\caption{Cleaning procedure  for the clan presented in Figure \ref{c07}. a)
$C=1$ b) $C=2$.}
\end{figure}

\newpage
\section{Studying the characteristics of the clan of ancestors through
  simulation results} \label{sec:simul}

\indent

We perform a 1,000 simulations for several values of $\lambda<
\lambda_c^\ast$. Conditioned on the event ``the point $(x,0)$ is
present at the free process'' (which has probability $1 -
\exp\{-\lambda \rho\}$), we observed the values of $N(\A^{(x,0)})$ --
total number of rectangles present in the clan.

 The expectations of this variable were estimated through the
 sample mean and compared then to the expected values for the
 associated branching process used to find the sub-criticality
 condition. 

The simulations were performed in two cases, when $\pi$ is the
$U(0,1)$ distribution and when $\pi$ is concentrated in one point
(fixed call length). From Garcia and Mari\'c (2003), the critical value for
$\lambda$ to assure sub-criticality is
\begin{itemize}
\item When $\pi=U(0,1)$
\begin{equation}
\lambda_c^\ast=\frac{1}{\frac{1}{2}+\sqrt{\frac{1}{3}}}~~\approx
0.9282
\end{equation}
\item When $\pi(d)=1$
\begin{equation}
\lambda_c^\ast=\frac{1}{2d}
\end{equation}
\end{itemize}

 Graphic1: N (uniform)
 \begin{figure}[p]\centering
 \psfig{file=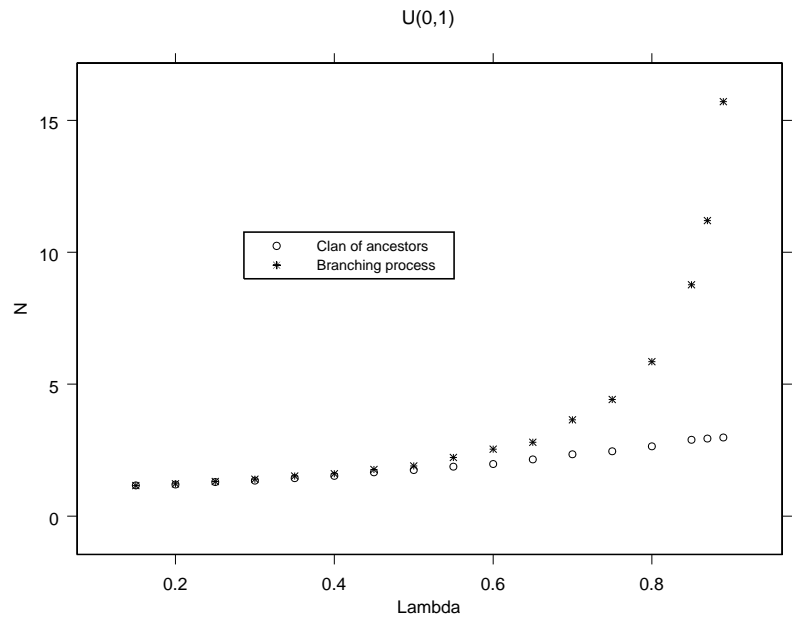,scale=.6}
 \psfig{file=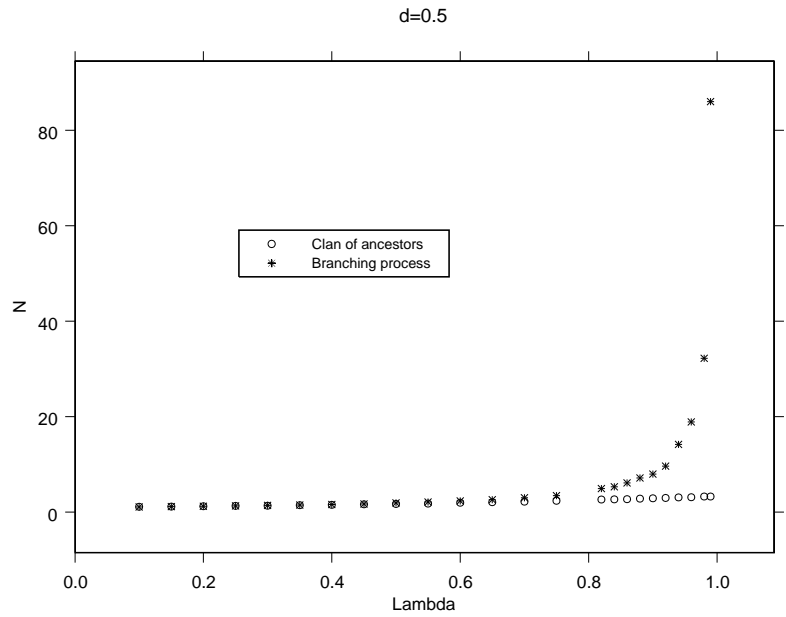,scale=.6}\label{Nlb}
 \caption{Expected total number of rectangles in the branching process
   and the clan of ancestors. a) $U(0,1)$ b) $d=0.5$}
 \end{figure}

Figure \ref{Nlb} shows that the branching process dominates the clan
of ancestors (we constructed then this way).  However, due to the fact
that in the branching process we can have subsequent generations of
rectangles to be born in the same area as the predecessor generations, 
  the number of
rectangles of the branching process grows much faster than the number of
rectangles of the clan of ancestors, as $\lambda$ increases. 

\subsection{Estimation of the critical value using simulations}

The purpose of this section is to study the behavior of the clan of
ancestors as $\lambda$ increases above $\lambda_c^\ast$. From Figure
\ref{Nlb} we can see that the critical
value obtained through the domination by a branching process
underestimates the true value of the finiteness of the clan of
ancestors. The idea behind these results is to generate samples for
increasing values of $\lambda$ and to study the total number of
rectangles. Our conjecture is that, close to the true critical value
$\lambda_c$ the total number of rectangles should grow exponentially
fast. Thus finding an assintote for $\E(N)$ would give us an estimate
of $\lambda_c$. This is true for the branching process, comparing the
value of $\E(N)$ as $\lambda$ approaches $\lambda_c^{\ast}$ in Figure
\ref{Nlb} we can see visually a vertical assintote at $\lambda_c^{\ast}$.

\begin{figure}[h,p]\centering

\psfig{file=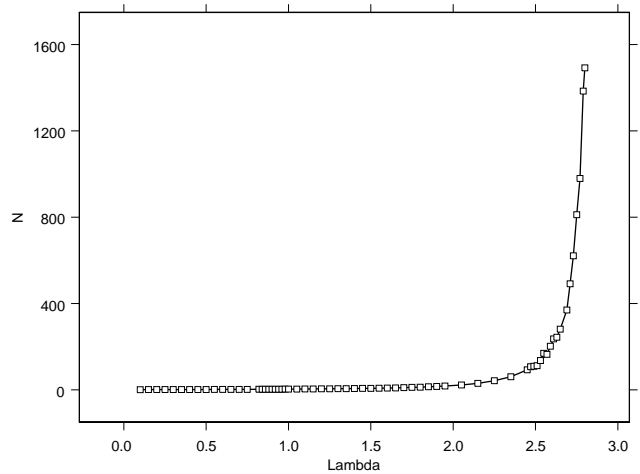,scale=.5}\label{N05}
\caption{ Expected number of rectangles in the clan of ancestors ($\E(N)$)
   for $d=0.5$} 
\end{figure}

From now on, for all distributions, we sampled 1,000 observations of
the clan of ancestors and computed $\bar{N}$ (the sample mean) for
each value of  $\lambda$. 
Figure \ref{N05} present the results for fixed length calls, $d=0.5$.

At first sight we see that there is an assintote close to  2.8. To be
more precise, we tried to find a root for
the equation  $1/\log(\bar{N})=0$ 
. We used a
degree 19 polynomial to approximate  $1/\log(\bar{N})$ and found a root
in $\lambda=2.8231$.

Several simulations were performed for several values of
$d$ just to get a more precise estimate for $\lambda_c$ since due to
the invariance of the Poisson process for fixed call length there is a
linear relationship among the critical values for all $d$, see Table
\ref{critconst}.

\begin{table}[h]

\begin{center}
\begin{tabular}{|c|ccccccccccc|} \hline
d  & 0.5 &  1.0  & 1.4  & 2.0 & 2.5 & 3.0 & 3.5 & 4.0 & 4.5 & 5.0 & \\ \hline

$\lambda_c$ & 2.8231 &  1.4193 &  1.0254&   0.7312 &  0.5682 &  0.4537
 & 0.3931 & 0.3530  &   0.3103 &  0.2833 &   \\ \hline 

\end{tabular}\label{critconst}
\end{center}
\caption{Critical value of $\lambda_c$ obtained through simulation for several call lengths}
\end{table}

Comparing the values of $\lambda_c^{\ast}=\frac{1}{2d}$ and
$\lambda_c$ we can see, as expected, a linear tendency and we can adjust a
regression model with no intercept using  least squares to get

\begin{equation}\label{hipotese}
\lambda_c = 2.8246\cdot \frac{1}{2d}.
\end{equation}

The question now is to perform the same comparison using different
random distributions for $\pi$. We simulated clan of ancestors for
several Beta distributions and compared  $\lambda_c^{\ast}$ and
$\lambda_c$. Table \ref{critrazne} presents these results along with
the ratio 
$\lambda_c / \lambda_c^{\ast}$. We can see that $\lambda_c \approx
2.82 \lambda_c^{\ast}$. Adjusting a least square model without any
intercept: 
\begin{equation}
\lambda_c= 2.8243\, \lambda_c^{\ast}.
\end{equation}

\begin{table}[h]

\begin{center}
 \begin{tabular}[t]{|c|c|c|} \hline
Distribution & \multicolumn{1}{|c|}{$\lambda_c$} & \multicolumn{1}{|c|}{$\lambda_c / \lambda_c^{\ast}$}\\ \hline U(0,1) & 2.6135 & 2.8157\\ \hline Beta(2,1) & 2.0888 & 2.8695 \\ \hline Beta(2,2) & 2.6746 & 2.8022
\\\hline
Beta(3,1) & 1.8597 & 2.8353\\ \hline Beta(3,2) & 2.3079 & 2.8444\\ \hline Beta(1,2) &3.7981 &2.8166\\ \hline
\end{tabular}\label{critrazne}
\end{center}
\caption{Critical value of $\lambda_c$ obtained through simulation for distributions $U(0,1)$ and Beta$(\alpha,\beta)$ }
\end{table}

\section{Dominating the clan of ancestors by a 2-generation branching
  process. Critical value.} \label{sec:ramificacao}
 \indent
 
 To construct the Galton-Watson process $b_n$ mentioned by the end of
 Section \ref{graphical}, Fern\'andez et al. (2001) used a multitype
 branching process $\B_n$, in the set of cylinders, which dominates
 $\A_n$.  To do this they looked ``backward in time'' and let
 ``ancestors'' play the role of ``branches''.  In particular, births
 in the original marked Poisson process correspond to disappearance of
 branches.  Like them, we reserve the words ``birth'' and ``death''
 for the original forward-time Poisson process. This construction can
 be done by enlarging the probability space and defining, for any
 given set $\{R_1,\dots,R_k\}$, \emph{independent} random sets
 $\B_1^{R_i}$ with the same marginal distribution as $\A_1^{R_i}$.
 The important point here is that
\begin{equation}
  \label{wt}
  \bigcup_{i=1}^{k}\A_1^{R_i} \subset \bigcup_{i=1}^{k}\B_1^{R_i}.
\end{equation}

The procedure defined by $\B_1$ naturally induces a multitype branching process
in the space of rectangles.  The $n$-th generation of the branching
process is defined by
\begin{equation}
  \label{an1}
  \B_n^R = \{ \B_1^{R'}: R'\in \B_{n-1}^R\}
\end{equation}
where for all $R'$, $\B_1^{R'}$ has the same distribution as
$\A_1^{R'}$ and are independent random sets depending only on $R'$. Then, for
all $n \geq 1$ it follows
\begin{equation}
  \label{wt1}
  \A_n^{R} \subset \B_n^{R}.
\end{equation}

 It is introduced then a multitype branching process in the set of calls $\G$.
For an initial rectangle $R$ with the basis of size $u$, $b_n^u$ is defined as the number of
rectangles in the n-th generation of $\B_n$ that have basis with length $v$:
\begin{equation}\label{eq:branching}
\b_n^{u}(v) = |\{ R^{'} \in \B_n^R |~ |\basis(R^{'})|= v\} |.
\end{equation}
The following relation is then established
\begin{equation}\label{relation}
     |\A_n^R| \leq |\B_n^R|= \sum_v b_n^u(v).
\end{equation}

\vskip5mm
 
We will have to further enlarge the set of types of $\B_n$ by creating a
branching process in the set of cylinders $\B^*_n$ with a new type
called $\color \in \{b ,g \}$ which is inherited by a cylinder not
only from the types of the first generation (parent) but also from the
2nd generation (grandparent).  The types $b, g$ are chosen to associate to
 \emph{black} and \emph{green}.\\
The idea is to identify some rectangles in the branching process $\B_n$
that could never be present in the $\A_n$. These rectangles we have in mind as ``black".
Then, if the total number of ``green" (not-black) rectangles is finite implies that the clan of
ancestors is also finite. We expect then that counting only ``green" rectangles one can
obtain  better critical value for $\lambda$. That is what we are going to prove in the following.\\

For a rectangle $V = (\xi_V, \tau_V, s_V, u_V)$ in the $n$th generation
of the branching process $\B_n$, we define an extra type to be its
color as follows. For $n=1$ or $2$, $\color(V) = g$. For $n > 2$, let
$R \in \B^*_{n-1}$ and $R' \in \B^*_{n-2}$ such that $V \in \B_1^{R}$ and $R \in
\B_1^{R'}$. Denote $R =  (\xi, \tau, s, u,c)$ and $R' =
(\xi', \tau', s', u',c')$, where $c = \color(R)$ and
$c'=\color(R')$. Therefore, 
\begin{itemize}
\item if $c = b$ then $\color(V) = b$, 
\item if $c = g$ then $\color(V) = b$ ~~if and only if~~ $(\xi_V, \tau_V +
  s_V) \in D(u_V, R, R')$ where $D(u_V, R, R') = L(\xi, u,
  \xi', u', u_V) \times [0,\tau']$ and 
\begin{equation}
\label{eq:a}
L(\xi, u,\xi', u', u_V) = [\max(\xi-u_V, \xi' - u_V), \min(\xi + u,
\xi' + u')]. 
\end{equation}
\end{itemize}
This way, for every rectangle is defined its bi-dimensional type
$\type(\cdot)\in \G \times \{b,g\}$.

Consider $d_n^{(w,c')}((u,c),(v,c''))$ to be a 2-generation multitype
branching process defined by
\begin{eqnarray}
  \label{eq:b}
  \lefteqn{d_n^{(w,c')}((u,c),(v,c'')) } \\
& &  = \, | \{V:V \in B_n^C \mbox{ and } R \in B_{n-1}^C;
  \type(C) = (w,c'), \type(R) = (u,c),\type(V) = (v,c'')\} |, \nonumber 
\end{eqnarray}
that is, it is the number of cylinders of type $(v,c'')$ in the $n$th
generation of a cylinder of type $(w,c')$ which have an ancestor in
the $n-1$ generation of type $(u,c)$.
 We can think of this branching
process as a first order branching process with enlarged state space given
by ${\cal S} = \{((i,u,c), (j,v,c''))\}$ and  offspring
distribution to be Poisson distributed.

From the definition of the process $d_n$ and the relation \reff{relation} it is clear that

\begin{equation}\label{eq:reldb}
     |\A_n^R| \leq \sum_{u,v} d_n^{(w,g)}((u,g),(v,g)) \leq \sum_v b_n^w(v).
\end{equation}

We suppose from now on that $\pi$ has discrete support.

\begin{propos} \label{prop:1}
In the process $d_n$ the mean number of offspring of type $(u,c'')$ of an individual of type
$(v,c)$ with the parent of type $(w,c')$ is given by
\begin{equation}\label{eq:meaneq}
m(((u',c'),(u,c));((u,c),(v,c''))) = \left\{ \begin{array}{ll}
                 0, & \mbox{ if }  c'=b,c=b,c''=g \mbox{ or } c'=g,c=b,c''=g  \nonumber \\
                     &\mbox{ or } c'=b,c=g,c''=b \mbox{ or } c'=b,c=g,c''=g\\
                 \lambda \pi(v)(u+v), &  \mbox{ if } c'=b,c=b,c''=b \mbox{ or }c'=g,c=b,c''=b  \nonumber 
                                    \end{array}
                                     \right.
\end{equation}
and
\begin{eqnarray} 
 m(((u',g),(u,g));((u,g),(v,b)) & \geq & \frac{1}{2} \lambda \pi(v)v,
                  \label{eq:mean} \\ 
 m(((u',g),(u,g));((u,g),(v,g)) & \leq & \lambda \pi(v)(u+\frac{v}{2}) \label{eq:mean2}
\end{eqnarray}
\end{propos}
\emph{Proof:}
For \reff{eq:meaneq} it is enough to observe that from the definition of Color  
follows that is impossible to have green children from black parents.
Since a black individual has all its children black, independently of its own parent type,
the number of its children has the same low as the first generation in branching process $b_n$
defined by \reff{eq:branching}. 
The total number of children of a green individual has also the same
low as the above one, but in this case it is possible to have children of both colors. Recall from \reff{media}
that $\E(b^{u}(v))= \lambda \pi(v)(u+v)$ so that
\begin{equation}
m(((u',g),(u,g));((u,g),(v,b)))+ m(((u',g),(u,g));((u,g),(v,g))) = \lambda \pi(v)(u+v).
\end{equation}
Therefore, to prove \reff{eq:mean} and \reff{eq:mean2} it is sufficient
to find the mean number of black children from green parent (and necessarily green grandparent): 
$m(((u',g),(u,g));((u,g),(v,b)))$. 

Consider two incompatible rectangles $R = (\xi, \tau, s, u,g) \in
\A_{2}$ and $R'=(\xi', \tau', s', u',g) \in \A_1$ such that
$R \in \A_1^{R'}$, that is $R$ is an ancestor of $R'$. Since $R'$ belongs
to the first generation, we have $\tau' < 0$ and $\tau'+s' >0$. By
construction of the clan of ancestors $\A_n$ and the branching process
$\B_n$, rectangles $R'$ in the first generation have  $\birth(R') < 0$
and $\death(R') > 0$. We
want to compute the $X(v, R,R')$- number of {\bf{black}} ancestors of $v$-type, namely the
number of those $V= (\xi_V,
\tau_V, s_V, v) \in \A_3$ with length $v$ such that $V \in \A_1^R$ and
$R \in \A_1^{R'}$.   Let $X_b(v, R,R')$ and $X_a(v, R,R')$ be the number
of such ancestors that died before and after time $t=0$, respectively.

As before, let  $D \equiv D(v, R, R') = L(\xi, u,
  \xi', u', v) \times [0,\tau']$ and 
\begin{equation}
L\equiv L(\xi, u, \xi', u', v) = [\max(\xi-u_V, \xi' - v), \min(\xi + u,
\xi' + u')]. 
\end{equation}
                              
Then, 
\begin{equation}
X(v,R,R') \,=\, X_b(v,R,R') + X_a(v,R,R')
\end{equation}

with
\begin{eqnarray}
\P( X_b(v,R,R')= k|u,\tau,\tau')= e^{-p_b \lambda
|D|}\frac{(-p_b\lambda |D|)^k}{k!}
\end{eqnarray}
and
\begin{equation}
\P(X_a(v,R,R')= k|u,\tau,\tau')= e^{-p_a\lambda
|L|}\frac{(-p_a\lambda |L|)^k}{k!}
\end{equation}
where $p_b$ is the  probability that a rectangle who died in the
area $D$ is really an ancestral of $R$ of type $v$, and
analogously for $p_a$. Therefore, 
\begin{eqnarray}
p_b&=&\pi(v)\P(Y-S<\tau) \\
&=& \pi(v) e^{\tau}(e^{-\tau'} -1)/{(-\tau')}
\end{eqnarray}
and
\begin{eqnarray}
p_a&=&\pi(v)\P(-S<\tau) \\
&=&\pi(v) e^{\tau}
\end{eqnarray}
where  $Y \sim U(0,\tau')$ (death time)   and  $S \sim \exp(1)$ (lifetime) are independent
random variables. We remind here that the rectangles are constructed
in the negative time (``past'') so $\tau, \tau' \leq 0$.  \\

Given $R$ and $R'$, $X_a(v,R,R')$ and $X_b(v,R,R')$ are independent random variables, thus $X(v,R,R')$ 
 is Poisson distributed  with mean
\begin{equation}\label{eq:poissmean}
\pi(v)(p_b \lambda |D|+ p_a \lambda |L|)= \pi(v) \lambda |L|
e^{(\tau-\tau')}.
\end{equation}

Notice that $\tau < \tau'$
and we are given the parent relation, so $\tau'-\tau$ has exponential
distribution with mean one. Hence, $\E(e^{-(\tau'-\tau)})= 1/2$.

Furthermore, since $\xi \in [\xi'-u,\xi'+u']$
\begin{equation}\label{eq:lgeq}
|L|= \min(\xi+u,\xi'+u') - \max(\xi-v,\xi'-v) ~~\geq~~ v
\end{equation}
From \reff{eq:poissmean} and the above inequality follows
\begin{equation}
 m(((u',g),(u,g));((u,g),(v,b)))= \E(X(v,R,R')) ~\geq ~\frac{1}{2}\lambda \pi(v) v 
\end{equation}

And  we have proved \reff{eq:mean} as desired. \hfill
$\square$

\paragraph{Sub-criticality:}

On account of the relation \reff{eq:reldb} we are interested in sub-criticality conditions for the green-type
population of the branching process $d_n$ described above.
Remember that green individuals may appear in the n-th generation  only as descendants of a
branch made of green individuals only. Therefore,
our aim is to establish conditions for the convergence of the series
 \begin{equation}\label{newseries}
\sum_{n\geq 1}\sum_{v_{n-1},u} m^{(n)}(((w,g),(v,g));((v_{n-1},g),(u,g)))
 \end{equation}
 where $m^{(1)}(\cdot\,;\cdot)=m(\cdot\,;\cdot)$ is given by \reff{eq:mean} and
 for $n > 1$ 
 \begin{eqnarray}
\lefteqn{ m^{(n)}(((w,c'),(v,c));((v_1,c_1),(u,c''))) } \\ &=&
\sum_{v_2,c_2}m^{(n-1)} (((w,c'),(v,c));((v_2,c_2),(v_1,c_1)))
m(((v_2,c_2),(v_1,c_1));((v_1,c_1),(u,c''))). \nonumber
\end{eqnarray}

To simplify notation let 
\begin{eqnarray}
\label{eq:ast}
m_\ast (v,u)& = \lambda \pi(u)(v+ \frac{u}{2})
\end{eqnarray}
and it follows from Proposition \reff{prop:1} that for all $u,v$ and independently of $w$
\begin{eqnarray}
m(((w,g),(v,g));((v,g),(u,g))) \leq m_\ast(v,u).
\end{eqnarray}

Therefore, using the notation above and \reff{media} we can easily estimate the series
\reff{newseries} to 
\begin{eqnarray}\label{seriebig}
\lefteqn{
\sum_{v_{n-1},u} m^{(n)}(((w,g),(v,g));((v_{n-1},g),(u,g)))\leq 
\sum_{u,{v_1}\ldots
{v_{n-1}}}m(v,v_1)m_\ast(v_1,v_2)\cdots m_\ast(v_{n-1},u) } \nonumber\\
&\quad= \lambda^n \sum_u\sum_{v_1}\cdots
\sum_{v_{n-1}}\pi(v_1)(v+v_1)\pi(v_2)\left( v_1+ \frac{v_2}{2}\right)
\cdots \pi(u)\left(v_{n-1}+ \frac{u}{2}\right)
\end{eqnarray}

 Observe that
\begin{eqnarray}
\sum_u \pi(u)\left(v_{n-1}+ \frac{u}{2}\right)= v_{n-1} +
\frac{\rho_1}{2} = f_1^\ast+ g_1^\ast v_{n-1}
\end{eqnarray}where $f_1^\ast= \frac{\rho_1}{2}$ and $g_1^\ast=1$,
and define inductively
\begin{eqnarray}
\sum_{v_{n-i+1}} \pi(v_{n-i+1})\left(v_{n-i}+
\frac{v_{n-i+1}}{2}\right)(f_{i-1}^\ast+ g_{i-1}^\ast v_{n-i+1})
 = f_i^\ast+ g_i^\ast v_{n-i}.
\end{eqnarray}

Then we have
\begin{equation}\label{matrix}
\left[
\begin{array}{rr}
f_{j+1}^\ast\\
g_{j+1}^\ast
\end{array} \right] =
\left[
\begin{array}{rr}
\rho_1 /2 & \rho_2/2 \\
1  & \rho_1
\end{array} \right]\cdot
\left[
\begin{array}{rr}
f_{j}^\ast\\
g_{j}^\ast
\end{array} \right]=
\left[
\begin{array}{rr}
\rho_1 / 2 & \rho_2/2 \\
1  & \rho_1
\end{array} \right]^{j} \cdot
\left[
\begin{array}{rr}
\rho_1/2\\
1
\end{array} \right]
\end{equation}
and consequently, the series \reff{seriebig} is equal to
\begin{eqnarray}\label{member}
\lambda^n \sum_{v_1} \pi(v_1)(v+v_1)
(g_{n-1}^\ast v_1 + f_{n-1}^\ast) = 
 \lambda^n
( v( g_{n-1}^\ast \rho_1 + f_{n-1}^\ast ) + \rho_2 g_{n-1}^\ast +
\rho_1 f_{n-1}^\ast).
\end{eqnarray}

\vskip8mm

In order to find $f_n^\ast, g_n^\ast$ we exponentiate $ T^\ast= \left[
\begin{array}{rr}
\rho_1 / 2 & \rho_2/2 \\
1  & \rho_1
\end{array} \right]
$. For this operation suffices the eigenvalues of $T^\ast$,
$\varepsilon_1$ and $ \varepsilon_2$  given by
\begin{equation}
\varepsilon_{1,2}=\frac{3\rho_1  \pm \sqrt{\rho_1^2 +
8\rho_2 }}{4}.
\end{equation}
and two corresponding normalized eigenvectors
\begin{eqnarray}
 \frac{1}{\sqrt{1+(\varepsilon_1-\rho_1)^2}}\left[
\begin{array}{c}
\varepsilon_1-\rho_1 \\
1
\end{array} \right],~~\frac{1}{\sqrt{1+(\varepsilon_2-\rho_1)^2}}\left[
\begin{array}{c}
\varepsilon_2-\rho_1 \\
1
\end{array} \right].
\end{eqnarray}
 From (\ref{matrix}) it follows

\begin{eqnarray}
f_n^\ast &=& \frac{1}{\varepsilon_1 -
\varepsilon_2}(\varepsilon_1^{n}(\varepsilon_1
- \rho_1) + \varepsilon_2^{n} (\rho_1-\varepsilon_2))\\
g_n^\ast &=& \frac{1}{\varepsilon_1 - \varepsilon_2} (
\varepsilon_1^n - \varepsilon_2^n).
\end{eqnarray}

\vskip7mm Using the fact that $|\varepsilon_2/\varepsilon_1|\leq
1$ in (\ref{member}) and the Cauchy-Hadamard formula, we obtain
that the radius of convergence of the series (\ref{seriebig}) is
\begin{equation}\label{bound1}
\lambda_c^{\ast\ast}= \frac{1}{\varepsilon_1}=\frac{4}{3\rho_1  + \sqrt{\rho_1^2 +
8\rho_2 }}
\end{equation}
 
Therefore, as long as $\lambda < \lambda_c^{\ast\ast}$ the series (\ref{newseries}) is absolutely convergent and 
consequently the green-type population of the process $d_n$ is almost surely finite. 

\vskip5mm Calculation being done for countable number of types is not
a real limitation.  Namely, let $V$ be the set of all possible types,
by the same argument as used in Proposition \ref{prop:1} the mean
number of green-type offspring in all generations  is less then
\begin{equation}
\sum_{n \geq 1}\int_V m^{(n)}_*(v,du)
\end{equation}
where
\begin{equation}
m^{(n)}_*(v,du) = \int_V m^{(n-1)}_*(v,dv')m_*(v',du)
\end{equation}
can be obtained inductively. 

Suppose that the distribution of the length of the calls is absolutely
continuous with respect to the Lebesgue measure and call $\pi$ its
density. We can write

\begin{equation}
 m_*(v,du)= \lambda \pi(du)(\frac{du}{2}+v) 
\end{equation}

and the computation is completely analogous to the discrete case.

\vskip8mm

  We summarize the results proved in this section in 
the following theorem:
\begin{thm}
If  $\lambda < \lambda_c^{\ast\ast}$ , where $\lambda_c^{\ast\ast}$ is given by (\ref{bound1}),
then  with probability one there is no backward oriented percolation in $\bf{R}$. 
\end{thm}

\vskip7mm
\noindent  \textit{Remark: If $\pi$ is the $U(0,1)$ density the condition becomes $\lambda <
1.247$}
\vskip3mm

Observe the following relations:
\begin{eqnarray}\label{1estim}
 \frac{2}{3\rho_1} \geq  \lambda_c^{\ast\ast} \geq 
 \frac{4}{3(\rho_1 + \sqrt{\rho_2})}= \frac{4}{3}
\lambda_c^\ast.
\end{eqnarray}
The upper estimate in (\ref{1estim})is achieved in the case of just one type (fixed call length).
The last inequality proves that a better bound for the critical
value is obtained.

\section{Conclusion}

This is one of the first works where the clan of ancestors algorithm
was implemented.  Berthelsen and M{\o}ller (2002) compared it to
the dominated CFTP introduced by Kendall and M{\o}ller (2000). Based
on simulation results, they show that the dominated CFTP is
better than the algorithm based on the clan of ancestors in the
particular case of a Strauss process
\begin{equation}
\label{eq:300ca}
 \mu_\Lambda(dN) = \frac{1}{Z_\Lambda} e^{\beta_1 N(\Lambda) + \beta_2
 S(N,\Lambda)}\, \mu^0_\Lambda(dN) 
\end{equation}
 defined on a unit square with $e^\beta_1 = 100$ and $e^\beta_2 = 0$
(the so-called hard-core process), $0.5$ and $1$ (a Poisson processes
with rate 100). This is obviously the case from the description of the
processes since the backward construction of BFA  stops when
the dominated Poisson process regenerates and usually the coupling of
CFTP is achieved before it in the finite case. However, it should be
noticed that the algorithm based on the clan of ancestors was designed
for sampling the infinite-volume process viewed in a finite window.
This seems to be a much more interesting and challenging problem which
has been studied in this work for the specific case of
one-dimensional loss networks with bounded calls. No comparison was
made to other perfect simulation schemes.

Moreover, we can see that the simulation of the invariant measure can
bring information about unknown variables related to the clan of
ancestors. The bound described in  Section \ref{sec:ramificacao} was found
by a better understanding of the simulation procedure. 

\paragraph {\bf Acknowledgments} We would like to thank Pablo Ferrari
for many fruitful discussions. This work  was partially funded by
FAPESP Grant 00/01375-8 and CNPq Grant 301054/93-2.

\end{document}